\documentclass{amsart} % in was article
\usepackage{amssymb}
\usepackage{amsmath}
\usepackage{amsfonts}

\begin{document}
%%%%%%%%%%%%%%%%%%%%%%%%%%%%%%%%%%%%%%%%%%%%%%%%%%%%%%%%%%%%%%%%%%%%%%
%	spaces for your own definitions follows
%%%%%%%%%%%%%%%%%%%%%%%%%%%%%%%%%%%%%%%%%%%%%%%%%%%%%%%%%%%%%%%%%%%%%%
\newtheorem{theo}{Theorem}[section]
\newtheorem{prop}[theo]{Proposition}
\newtheorem{lemma}[theo]{Lemma}
\newtheorem{exam}[theo]{Example}
\newtheorem{coro}[theo]{Corollary}
\theoremstyle{definition}
\newtheorem{defi}[theo]{Definition}
\newtheorem{rem}[theo]{Remark}

%\renewcommand{\theequation}{\mbox{\arabic{section}.\arabic{equation}}}

%letters - added these
\newcommand{\Bb}{{\bf B}}
\newcommand{\Nb}{{\bf N}}
\newcommand{\Qb}{{\bf Q}}
\newcommand{\Rb}{{\bf R}}
\newcommand{\Zb}{{\bf Z}}
\newcommand{\Ac}{{\mathcal A}}
\newcommand{\Bc}{{\mathcal B}}
\newcommand{\Cc}{{\mathcal C}}
\newcommand{\Dc}{{\mathcal D}}
\newcommand{\Fc}{{\mathcal F}}
\newcommand{\Ic}{{\mathcal I}}
\newcommand{\Jc}{{\mathcal J}}
\newcommand{\Lc}{{\mathcal L}}
\newcommand{\Oc}{{\mathcal O}}
\newcommand{\Pc}{{\mathcal P}}
\newcommand{\Sc}{{\mathcal S}}
\newcommand{\Tc}{{\mathcal T}}
\newcommand{\Uc}{{\mathcal U}}
\newcommand{\Vc}{{\mathcal V}}

\newcommand{\ax}{{\rm ax}}
\newcommand{\Acc}{{\rm Acc}}
\newcommand{\Act}{{\rm Act}}
\newcommand{\ded}{{\rm ded}}
\newcommand{\Gm}{{$\Gamma_0$}}
\newcommand{\ID}{{${\rm ID}_1^i(\Oc)$}}
\newcommand{\PAP}{{${\rm PA}(P)$}}
\newcommand{\ACA}{{${\rm ACA}^i$}}
\newcommand{\RefP}{{${\rm Ref}^*({\rm PA}(P))$}}
\newcommand{\RefS}{{${\rm Ref}^*({\rm S}(P))$}}
\newcommand{\Rfn}{{\rm Rfn}}
\newcommand{\tar}{{\rm Tarski}}
\newcommand{\UNFA}{{${\mathcal U}({\rm NFA})$}}

\author{Nik Weaver}

\title [The concept of a set]
       {The concept of a set}

\address {Department of Mathematics\\
          Washington University in Saint Louis\\
          Saint Louis, MO 63130}

\email {nweaver@math.wustl.edu}

\date{\em December 15, 2006}

\begin{abstract}
Metaphysical interpretations of set theory are either inconsistent or
incoherent. The uses of sets in mathematics actually involve three
distinct kinds of collections (surveyable, definite, and heuristic),
which are governed by three
different kinds of logic (classical, intuitionistic, and minimal). A
foundational system incorporating this analysis and based on the
principles of mathematical conceptualism \cite{W1} accords better
with actual mathematical practice than Zermelo-Fraenkel set theory does.
\end{abstract}

\maketitle

%%%%%%%%%%%%%%%%%%%%%%%%%%%%%%%%%%%%%%%%%%%%%%%%%%%%%%%%%%%%%%%%%%%%%%%
%	Please insert the article body now
%%%%%%%%%%%%%%%%%%%%%%%%%%%%%%%%%%%%%%%%%%%%%%%%%%%%%%%%%%%%%%%%%%%%%%%

\section{Metaphysical interpretations}

\subsection{Sets in ordinary language}

Elementary introductions to set theory tend to give the impression that the
concept of a set is trivial, something with which we are already thoroughly
familiar from everyday life. We may be told that such things as a flock of
birds, a deck of cards, or a pair of apples are examples of sets.

This immediately seems strange because sets in the mathematical sense
are supposed to be abstract objects not existing in space and time,
whereas it is hard to believe that a simple assertion about, say, a flock
of birds carries any significant metaphysical content.

Slater (\cite{Sla}, Section II) analyzed the way collective expressions
are actually used in ordinary language and showed that they indeed
express no metaphysical content. A deck of cards exists in
space and time; it is a physical object composed of cards in the
same way that a house is made out of bricks, something philosophers
call a ``mereological sum''. It is not a {\it set} of cards in the
mathematical sense, with cards as elements. The same is true of
such things as a bunch of grapes, a herd of cows, etc. (You don't eat
{\it the elements of} a bunch of grapes, you eat the bunch of grapes. And
``the fact that such collections are mereological sums is also shown
by the fact that shoals, herds, packs, tribes, and the like, are
located and can move around in physical space, just like their
members'' (\cite{Sla}, p.\ 61).)

In contrast, the word ``pair'' in the phrase ``pair of apples'' is a
numerical measure word; to say that you ate a pair of apples merely
tells how many apples you ate. It is analogous to the word
``half'' in the phrase ``half a loaf of bread''. And ``in `There is
a half of a loaf' there is obviously no reference to something other
than bread; there is not, in addition, reference to one of a range
of mysterious, further objective entities, `halves', `quarters',
`parts', etc. There is merely a specification of how much of a
loaf there is'' (\cite{Sla}, p. 62).

The idea that the mathematical concept of a set is obvious
and in no need of any special explanation
is not correct. This is sometimes noted in the philosophical
literature on mathematical foundations, but it seems the conclusion
that is usually drawn is not that the set concept is bankrupt but
rather that it merely needs some more sophisticated interpretation.

\subsection{Inconsistency versus incoherence}

But however we might interpret set language, it is clear that a
straightforward belief in the actual existence of a well-defined
objective world of sets immediately gives rise to the classical
paradoxes of naive set theory. Thus, naive metaphysical
interpretations of set language are inconsistent. Is there any cogent
alternative metaphysical interpretation which evades the paradoxes?

The standard answer is that this is to be done in terms of the {\it iterative
conception} of sets, according to which sets are not to be thought of as
arbitrary collections, but rather as being hierarchically built up from
the empty set via the two operations of (1) forming a subset of a given
set and (2) collecting together all subsets of a given set, i.e.,
forming its power set. This suggestion can be interpreted in two ways.
On one reading, it marks a distinction between two concepts, ``set''
and ``collection'', with sets effectively being just those collections
which appear in the cumulative hierarchy. But this is nothing
more than a change in terminology, so that the classical paradoxes about
sets are not defused but simply become
paradoxes about collections. Although the paradoxes no longer directly
invalidate set theory, they now do so indirectly by showing that the
notion of a collection, in terms of which sets are defined, is itself
inconsistent. Thus this interpretation of the set concept is incoherent.

The iterative conception can also be understood in a different way, not as
differentiating sets out of a background universe of collections, but rather
as clarifying the concept of a collection as something which must be in some
sense ``formed'' out of elements that in some sense exist ``before'' it does.
Exactly what this means for abstract objects not existing in space and time
is hard to pin down, so its status as a {\it clarification} is questionable,
to say the least. The point is apparently that the existence of a set has
something to do with its being constructible, in some obscure sense.
What makes this interpretation really incoherent is the fact that each stage
of the process by which sets are supposed to be built up involves forming a
power set, which in the case of infinite sets is an absolutely non-constructive
operation. In other words, the set-theoretic universe is thought of as being
{\it built up} in an iterative process, yet we pass from each stage of this
process to the next in a completely non-constructive way. The power set of the
preceding level is not constructed in any sense whatever, it simply appears.

The source of this incoherence is obviously the fact that axiomatic set
theory was developed not in a philosophically principled way, but rather
in an opportunistic attempt to preserve as much naive set theory
as possible without admitting any obvious inconsistencies. This was
accomplished by, in effect, adding just enough constructivity to the
naive picture to avoid the standard paradoxes. The result is a
nonsensical jumble of constructive and non-constructive ideas.

If the temporal (``in stages'') and constructive metaphors, which are
highly dubious anyway in light of the purported abstract nature of sets,
are left out of the iterative conception, then it loses its force as a
resolution of the paradoxes. One is simply left with a bald assertion
that sets are layered in a hierarchy with no explanation as to why
this is so. And we still have the problem of justifying the power set
axiom.

It seems that the fundamental error in all metaphysical interpretations
of set theory is the reification of a collection as a {\it separate
object}, and that this is done as a result of a series of grammatical
confusions \cite{Sla}. This reification is the ultimate source of the
paradoxes, and once one accepts it there is no cogent way to avoid them.

\subsection{The philosophical dilemma}

The basic dilemma faced by philosophers of mathematics
was summarized well by Hazen: ``On the one hand, the notion of a set is
central to modern mathematics $\ldots$ On the other hand, we seem unable
to obtain for set theory the kind of metaphysical or epistemological
legitimacy that would come from a characterization of the sort of entity
a set is or an account of how we become acquainted with them''
(\cite{Haz}, p.\ 173).

Because the task of legitimizing axiomatic set theory has seemed so
crucial, a great deal of effort has been expended towards this goal.
Various authors have either proposed to defend set theory on
its own (metaphysical) terms, or tried to justify it in
some weaker sense by more roundabout methods. However, it seems fair
to say that none of these attempts has met with general approval.

But the dilemma rests on a faulty assumption. Although the {\it language} of
set theory is used throughout mathematics in an elementary way, the actual
{\it discipline} of axiomatic set theory is not central to modern mathematics.
In fact it is quite peripheral, making only occasional and relatively minor
contact with mainstream areas. Indeed, virtually all modern mathematics
outside set theory itself can be carried out in formal systems which are
far weaker than Zermelo-Fraenkel set theory and which can be justified in
very concrete terms without invoking any supernatural universe of sets (see
\cite{W2} or \cite{W5}, and also \cite{W3} for borderline
cases). Thus, axiomatic set
theory is not indispensable to mathematical practice, as most philosophers
of mathematics have apparently assumed it to be. It is one arena in which
mathematics can be formalized, but it is not the only one, nor even
necessarily the best one (see Section 3.3 below).

The point of view taken here could be expressed by moving the word
``seem'' in the passage quoted above: we {\it are} unable to obtain
metaphysical or epistemological legitimacy for set theory, but it only
{\it seems} that sets are central to modern mathematics. I discuss this
issue further in \cite{W4}.

\section{Collections}

\subsection{Proxies for sets}

We must reject as nonsensical the idea that sets literally exist as some sort
of mysterious non-physical entities. But this does not mean that we cannot
make sense of any kind of set language. For example, we do not need to
be set-theoretic platonists to grasp the meaning of the statement that every
nonempty set of natural numbers has a least element. The same idea
could be expressed by saying that if either a 0 or a 1 appears
in each cell of a one-way infinite tape, and at least one cell contains 1,
then there is a first cell that contains 1. This formulation is clumsier than
the first one and it introduces extraneous notions (the symbols 0 and 1, the
image of an infinite tape), but it does show us that we can understand what
is expressed by the first statement in a way that does not require us to
think of sets of natural numbers as actual objects.

This is somewhat analogous to the difference between ``two plus three is five''
and ``two apples plus three apples is five apples''. Here too the second
formulation is clumsier and unnecessarily specific. But again, it shows that
we do not need to assume the literal existence of a platonic world
of numbers in order to make sense of basic arithmetic. In other words, we do
not gain any substantive mathematical content by supposing that number words
literally refer to some sort of abstract objects.

Thus, the idea is that we can find meaning in at least some language
involving fictional ``abstract objects'' like numbers and sets by using
actual physical objects (apples) or possible physical objects (infinite
tapes) as proxies for them. Of course philosophical questions can be raised
about the idea of an infinite tape, but the point is that we avoid the far
more serious difficulties attaching to set theory. The naive notion of an
infinite tape is not obviously paradoxical, and Hazen's questions (What
sort of entities are they? How do we become acquainted with them?)
hardly have the same force they have against abstract sets.

The moderately substantial set-theoretic system used in \cite{W2} could be
interpreted in a similar way, in this case using formal expressions as
proxies for sets. In some ways this example is even less problematic because
the structure in question,
$J_2$, is countable and the formal expressions are finite. $J_2$ can be
thought of as a toy model for the Cantorian universe, and despite its
small size it is, surprisingly, rich enough to permit the development of
ordinary mainstream mathematics, even up to subjects like measure theory
and functional analysis. This was shown in some detail in \cite{W2}.

Because it is so concrete, if we agree to practice mathematics only within
$J_2$ then philosophical concerns about sets largely disappear. However,
although $J_2$ is an attractive structure and is quite adequate for
ordinary mathematics, we can easily extend its construction to obtain
richer toy universes. Reasoning about $J_\alpha$ for arbitrary $\alpha$
then reintroduces philosophical issues because we cannot model the class
of ordinals with a single physically instantiated structure. This forces
us to come to grips with the general notion of a collection outside of any
particular system of proxies.

\subsection{Types of collections}

There are three essentially different kinds of collections. First, we have
collections that could be instantiated in some possible world by a physical
structure in which the elements of the collection appear as discrete
components. I call this kind of collection {\it surveyable}. In principle
we could exhaustively search through all the individuals in such a
collection.

What counts as a possible physical structure is obviously open to debate.
According to finitism it is precisely the finite collections that are
surveyable. In classical set theory, on the other hand, a collection is
surveyable if and only if it is a set. It might take a transfinite amount
of space to contain the elements of a set, or a transfinite amount of time
to survey them, but this would be seen as falling within the realm of
logical possibility.

The next kind of collection is {\it definite}. This is a much broader
category. We no longer require that we be able to search through all the
individuals in the collection; instead, we ask only that the statement that
an individual belongs to the collection have a fixed well-defined meaning.
For example, the prime numbers are finitistically definite, since we can
finitistically test any natural number for primality. Thus the statement
that a number is prime has, finitistically, a fixed well-defined meaning.
In classical set theory ``definite collection'' is synonymous with ``class''.

Definite collections generally cannot be physically modelled in the way that
surveyable collections can. So the most natural proxies for definite
collections will be the predicates that define them.

We do not assume that membership in a definite collection is decidable.
For example,
consider the values of $k$ such that every graph with chromatic number
$k$ has a $k$-clique minor. On the face of it we cannot finitistically
test whether a given natural number has the stated property, because
determining the truth of this condition for any given value of $k$
apparently requires us to quantify over infinitely many graphs. But the
condition still has a fixed well-defined meaning for each value of $k$.
As we build up a successively larger repertoire of natural numbers our
understanding of what the condition means does not change. Thus, this
collection is finitistically definite but (on its face) not decidable.
In classical set theory we are in the same position with regard to classes
that are defined by means of conditions that involve quantification
over the entire set-theoretic universe (or over any proper class).

Our notion of definite collections may be contrasted with Dummett's notion of
``indefinitely extensible'' concepts \cite{Dum}. One trivial difference is
that we include surveyable collections among definite collections, whereas
Dummett's indefinitely extensible concepts are supposed to be
complementary to definite concepts. But any collection that is, in the
above sense, definite but not surveyable would have to have the property
that it is not exhausted by any surveyable subcollection. That is, given
any surveyable collection that is contained in a definite collection
which is not surveyable, there would have to be a new individual that
belongs to the latter but not the former. This looks similar to
Dummett's condition that ``if we can form a definite conception of a
totality all of whose members fall under the concept, we can, by
reference to that totality, characterize a larger totality all of
whose members fall under it'' (\cite{Dum}, p.\ 441). However, 
Dummett seems to require that there be an explicit prescription
for producing these new elements and we do not assume this. Also,
Dummett apparently takes all concepts to be decidable, which, as I
just indicated, I do not do.

In the case of a collection that is definite but not surveyable, we
can go beyond any surveyable subcollection, but the collection as a
whole is still well-defined. The last kind of collection, which I call
{\it heuristic}, goes one step further: any definite subcollection can
be extended. For example, the collection of all well-defined predicates
is heuristic because the property of being a well-defined predicate
does not itself have a fixed well-defined meaning. To see this let $P$
be any well-defined predicate with the property that $P(x)$ implies that
$x$ is a well-defined predicate. Then let $Q$ be the predicate ``$P(x)$
and $x$ does not hold of itself''. Since $P$ is well-defined so is $Q$,
but if $P$ holds of $Q$ then we immediately get a contradiction: if $Q$
holds of itself then it does not hold of itself, and vice versa. We
cannot imagine a state of affairs in which the assertion that $Q$ holds
of itself has a truth value, so $Q$ cannot be a well-defined predicate,
contradicting the assumption that $P$ holds of $Q$. So $P$
does not exhaust all well-defined predicates. This shows that any definite
collection of well-defined predicates can be extended. Other basic
examples of heuristic collections include the collection of all valid
proofs and the collection of all valid definitions.

By definition, we cannot give a precise meaning to the notion of
a heuristic collection. But we would still like to be able to reason
about --- in particular, to quantify over --- all well-defined predicates,
or all valid proofs, or all valid definitions. (For instance, we want to
be able to say: if there is a valid proof of $A$ and a valid proof of
$B$ then there is a valid proof of $A \wedge B$.) The classical semantic
paradoxes might be taken to show that this is an unreasonable wish.
But in fact these paradoxes only arise when heuristic concepts are
treated as if they were definite. Once the distinction between definite
and heuristic concepts is recognized the semantic paradoxes are easily
resolved \cite{W6}.

I will say that a variable is surveyable, definite, or heuristic
according to whether it is supposed to range over the individuals
belonging to a surveyable, definite, or heuristic collection.

In Section 3 I will make a case that there is a reasonable version of the
power set operation for surveyable and definite collections, but that the
result of this operation when it is applied to a surveyable collection
in general is merely definite, and the result when it is applied to a
definite collection in general is merely heuristic. There is no meaningful
version of the power set operation for heuristic collections.

\subsection{Systems of logic}
Legitimate forms of logical reasoning vary depending on whether the
language one is using expresses surveyable, definite, or
heuristic properties. I will argue that the appropriate
corresponding forms of logical reasoning are respectively classical,
intuitionistic, and minimal. Similar suggestions have been made before;
most notably, Dummett has proposed that intuitionistic logic should be
used when reasoning about the individuals falling under an indefinitely
extensible concept.

The differences between the three logics are most elegantly expressed
in terms of natural deduction \cite{Pra}. To present the logical rules
as simply as possible we introduce a symbol $\perp$ for ``falsehood'' and
regard $\neg A$ as an abbreviation of $A \to \perp$. Informally, the rules
for classical logic are:
\begin{enumerate}
\item
Given $A$ and $B$ deduce $A \wedge B$; given $A \wedge B$ deduce $A$ and
$B$.
\item
Given either $A$ or $B$ deduce $A \vee B$; given $A \vee B$, a proof of
$C$ from $A$, and a proof of $C$ from $B$, deduce $C$.
\item
Given a proof of $B$ from $A$ deduce $A \to B$; given $A$ and $A \to B$
deduce $B$.
\item
Given $A(x)$ deduce $(\forall x)A(x)$; if the term $t$ is free for $x$, given
$(\forall x)A(x)$ deduce $A(t)$.
\item
If the term $t$ is free for $x$, given $A(t)$ deduce $(\exists x)A(x)$; if
$y$ does not occur freely in $B$, given $(\exists x)A(x)$ and a proof
of $B$ from $A(y)$ deduce $B$.
\item
Given $(A \to \perp) \to \perp$ (i.e., $\neg\neg A$) deduce $A$.
\end{enumerate}
In intuitionistic logic the final rule is weakened to the ex falso law
``given $\perp$ deduce $A$'' (for any formula $A$), and in minimal logic
it is dropped altogether. (For a more precise exposition of this system
of reasoning see Chapter 1 of \cite{Pra}.)

\subsection{Classical versus constructive truth}
Formal assertions can be interpreted either classically or constructively.
The classical interpretation is formulated in terms of an assumed underlying
reality: we understand the statement that $A$ is true just to mean that what
$A$ asserts {\it is the case}. Thus, for example, we interpret $A \to B$
to mean that any evaluation of the free variables which makes $A$ the
case also makes $B$ the case. In contrast, the constructive interpretation
of formal assertions has to do with {\it what can be proven}. Under this
interpretation, to say that $A$ is true is to say that we are entitled to
assert $A$, i.e., that we can prove $A$. For example, here we are entitled
to assert $A \to B$ --- that is, we can prove $A \to B$ --- precisely if we
have a construction that will convert any proof of $A$ into a proof of $B$.

Although the classical interpretation is arguably more straightforward,
I claim that it cannot be used in the heuristic setting where definitions
are not fixed and hence there is no well-defined underlying reality; that in
the surveyable setting, on the other hand, there is no effective difference
between classical and constructive truth; and thus that the only case where
a real choice has to be made is in the definite setting.

First, for statements that involve only surveyable variables and predicates,
there is operationally no difference between the classical and constructive
interpretations. In principle, we could determine
the truth value of any such statement simply by performing an exhaustive
search (if there are multiple quantifiers, a finite series of exhaustive
searches). In other words, any sentence that is classically true could in
principle be proven and hence would also be constructively true. So in this
kind of setting classical logic is appropriate. In particular, the law of
excluded middle (the statement $A \vee \neg A$, for all formulas $A$, which
is generally absent from intuitionistic logic) is constructively valid here.

For statements which quantify over variables that are definite but not
surveyable, the situation is more subtle. We are still entitled to assume
an underlying notion of {\it being the case} which would support a classical
interpretation of truth. However, in general we should not expect that every
classically true sentence will be provable. That is, we may, even in principle,
have no way to determine the truth value of a statement that quantifies over
definite variables. Consequently, either classical or intuitionistic logic
can be used. But adopting a classical interpretation in a definite but not
surveyable context does not seem very helpful. All we gain is the existence
of truth values that in principle could never be known, and this cannot have
any substantive consequences.

The fact that classical logic is an option does, however, imply that under
the constructive interpretation of truth the ex falso law is
valid. Even if the classical truth values of some statements might not be
knowable, we can still incorporate their existence into our notion of a
valid proof. That is, we can insist that a proof cannot be valid if it yields
a conclusion that is classically false. Given this restriction, even though
we cannot access all classical truth values, we can still affirm
that $\bot$ is not provable. Thus the ex falso law $\bot \to A$ is
constructively valid because it is vacuously the case that we can convert
any proof of $\bot$ into a proof of $A$.

The preceding argument no longer holds in the heuristic case. In this
setting it does not make sense to suppose an underlying classical reality.
(If we could do this, we would be dealing with definite concepts.) So we
have to interpret truth constructively. Furthermore, we have no sharp global
distinction between valid and invalid reasoning. Thus we can no longer argue
that $\bot$ has no proof, and hence that the ex falso law is valid, because
this
assumes that constructive reasoning is sound (as otherwise a false statement
might indeed have a proof). Now one may find the validity of constructive
reasoning perfectly compelling, but the point is that this must be
decided before the justification of ex falso can proceed. The problem
is that this leads to an infinite regress if one is arguing that ex falso
itself ought to be included among the axioms of constructive reasoning.

This difficulty is genuine in the heuristic setting where we are not
dealing with fixed concepts with well-defined meanings. Here any
axioms we adopt will to some extent play a definitional role, in contrast
to the definite case where we need only to choose axioms which
accurately embody concepts that are already well-defined. In other words,
formal systems modelling definite concepts are merely descriptive, whereas
systems modelling heuristic concepts may be at least partly prescriptive.
In the latter context we lack the guarantee of consistency that comes from
correctly modelling an already well-defined concept, and there is a real
danger that adopting ex falso could introduce an inconsistency.

To summarize: if we adopt a classical interpretation of truth then we
can use classical logic in both the surveyable and definite settings, but
we cannot reason about heuristic collections. If we adopt a constructive
interpretation of truth then we must use classical logic in the surveyable
case, intuitionistic logic in the definite case, and minimal logic in the
heuristic case.

\section{Conceptualist mathematics}

\subsection{Supertasks}

The basic tenets of mathematical conceptualism \cite{W1} can be
summarized as follows: mathematics is the study of logical possibility;
conceivability is a sufficient condition for logical possibility; and
constructions of length $\omega$ are conceivable.

The first point can be contrasted with the view that mathematics is about
{\it physical} possibility. This view is expressed, for example, when
philosophical conclusions about mathematics are drawn from purported
upper bounds on the number of elementary particles in the universe. The
conceptualist objection to arguments like these
would be that any computational constraints we
face as a result of special features of our universe are not relevant
to the concept of mathematical truth. They may be interesting and
important in many other ways but they have no bearing on this issue.

Some sort of physicalist attitude seems to underlie constructivist ideas
generally, in particular the idea that only those constructions which
terminate after a finite number of steps are legitimate because it is
only these which could actually be performed.
This suggestion draws a very clean line between what is
acceptable and what is not, but it is suspect for that very reason. What
``could actually be performed'' is not a sharp concept, and putting all
finite computations into this category already involves a substantial
idealization. This leads to the question why we should not allow the
further idealization to computations of length $\omega$.

Various commentators have felt that computations of length $\omega$ ---
``supertasks'' --- clearly fall within the realm of logical possibility
\cite{Moo}.
If this point is granted, then in the terminology of Section 2 we should
conclude that the natural numbers are surveyable. Moreover, the truth value
of any statement of first order number theory could be mechanically decided
using infinite truth tables by a computation
of length $n\cdot \omega$ for some $n$ (i.e., a finite series of supertasks),
so that classical logic is indeed appropriate to this setting, in line with
the suggestion made in Section 2.4.

More generally, any countable structure familiar from ordinary mathematics
ought to be surveyable on the preceding analysis. Of course we need to
specify how such a structure is to be physically modelled, but in worlds
that permit supertasks there should be no serious difficulty in finding
appropriate models. At this level logic is classical and abstract set
theory is wholly absent.

\subsection{Uncountable sets and proper classes}

Thus, at the countable level conceptualist mathematics is practically
indistinguishable from classical mathematics. We are ``platonists about
number theory''. However, this is no longer true at the level of (what
are classically thought of as) uncountable sets. We cannot regard
uncountable structures, even familiar ones like the real line, the
power set of $\omega$, or $\aleph_1$, to be surveyable. This is because
we have no clear conception of what it would be like to look through all
real numbers, or all countable ordinals, even if supertasks were allowed.

The surveyability of collections like these is widely taken to be
self-evident. But this is just a consequence of failing to draw a clear
distinction between surveyability and definiteness. If one lacks this
conceptual distinction then one is bound to conclude that $\Pc(\omega)$
is {\it the same kind of thing} as $\omega$ since, after all, both are
collections. This line of thought leads straight into the bizarre
theology of Cantorian set theory, with its remote cardinals that have
no connection to mainstream mathematics.

In fact the surveyability of $\Pc(\omega)$ is an extraordinary claim,
completely different in character from the idea that $\omega$ is surveyable.
It is inherent in our conception of natural numbers that we can imagine
searching through them one at a time.  But we have no clear intuition
for how one would go about exhaustively searching through all sets of
natural numbers. It is only to be expected that mathematics that
requires the assumption that $\Pc(\omega)$ is surveyable would be highly
pathological, have no meaningful scientific applications, and make little
connection to mainstream mathematics generally.

The point can be made in terms of possible worlds. Can we clearly conceive
of a world in which (representatives of) all real numbers are present and
surveyable? The L\"owenheim-Skolem theorem suggests that we cannot. Any
description we could give of such a world as an abstract structure would
already be satisfied by some countable substructure, which evidently would
not contain all real numbers.

However, classically uncountable sets like $\Pc(\omega)$ and $\aleph_1$
are at least conceptualistically definite. If we accept that $\omega$ is
surveyable then we should also accept that any infinite sequence of $0$'s
and $1$'s is surveyable, leading to the conclusion that the predicate ``is
an infinite sequence of $0$'s and $1$'s'' is in principle decidably true or
false. This shows that the power set of $\omega$ can be conceptualistically
regarded as a definite collection. The real line can be handled in essentially
the same way, either identifying real numbers with their binary expansions,
or in terms of Dedekind cuts.

Once we accept that the real line is definite, it is a short step to the
same conclusion for all of the standard spaces that one encounters in
ordinary mathematics. In
some cases this will require fairly straightforward coding. For instance,
the real Banach space $C[0,1]$ may be regarded as definite because
elements of this space can be identified with uniformly continuous
functions from $[0,1] \cap \Qb$ to $\Rb$; such a function is specified
by a countable family of real numbers, which is easily seen to be
definite, or it can be encoded as a single real number, allowing a
reduction to the fact that $\Rb$ is definite. For details see \cite{W5}.

In general, any classical set whose elements can be reasonably
encoded as real numbers will be conceptualistically
definite. For instance, the set of all separable manifolds and the
set of all separable Banach spaces fall in this category via encodings
familiar from reverse mathematics \cite{Sim}.

Broadly speaking, the proper classes that appear in ordinary mathematics
(e.g., the class of all manifolds, or the class of all Banach spaces)
are conceptualistically realized as heuristic collections.

\subsection{Comparison with Zermelo-Fraenkel set theory}

In classical set theory every set has a power set. Conceptualistically, the
situation is not so simple. The example of $\omega$ and $\Pc(\omega)$
illustrates the principle that for every surveyable collection there is a
definite collection which plays the role of its power set. ``Subsets'' of
$\omega$ can be represented by sequences of $0$'s and $1$'s, which are
surveyable just like $\omega$, so that the collection of infinite sequences
of $0$'s and $1$'s is definite. But it is not surveyable.

The point is that sequences of $0$'s and $1$'s play the role of proxies
for surveyable subcollections of $\omega$. In general, given any surveyable
collection we should be able to set up a natural system of proxies for
surveyable subcollections, such that the collection of all surveyable
subcollections is definite. {\it The surveyable subcollections of any
surveyable collection constitute a definite collection.}

(We may note here that there exist non-surveyable but definite subcollections
of $\omega$, for instance the Church-Kleene ordinal notations, as well as
non-definite but heuristic subcollections of $\omega$, for instance the
G\"odel numbers of the well-defined predicates.)

Similarly, {\it the definite subcollections of any definite collection
constitute a heuristic collection}. Since the notion of a well-defined
predicate is merely heuristic, so is the notion of a definite subcollection
of a given collection, in general.

The process stops here. There is no meaningful sense in which one can talk
about arbitrary heuristic subpredicates of a given heuristic predicate,
because the general concept of a heuristic predicate is not even well-defined.

The accord between conceptualist philosophy and
mainstream mathematics is striking. In general, discrete
objects in ordinary mathematics are conceptualistically represented as
surveyable; continuous objects are represented as definite; and proper
classes are represented as heuristic. In Zermelo-Fraenkel set theory the sets
$\omega$, $\Pc(\omega)$, and $\Pc(\Pc(\omega))$ are just the first three
terms of an infinite sequence. But {\it all} mainstream mathematics takes
place, or can be construed as taking place, in these sets. To fully appreciate
this point, consider that classically the power set operation can be
iterated not just any finite number of times, but any transfinite number
of times (taking unions at limit stages) --- even uncountably many times.
Even, say, a measurable cardinal number of times, if one believes in measurable
cardinals. Thus ordinary mathematics is {\it grossly} discordant with the
Cantorian universe, while it fits the conceptualist picture perfectly.

The conceptualist universe is much smaller than the Cantorian universe,
but what is lost is a vast realm of set-theoretic pathology that plays no
essential role in mainstream mathematics. This is particularly seen in
fields where the basic definitions allow objects of arbitrary cardinality.
For instance, Banach spaces can be arbitrarily large. Yet nearly all
Banach spaces of real mainstream interest are either separable or the duals
of separable spaces, and hence fall into the conceptualist framework. There
are a few important spaces that are neither separable nor the duals of
separable spaces (such as spaces of almost periodic functions, or CCR
algebras, or the Calkin algebra), but these rare exceptions invariably
still turn out to lie within the conceptualist framework. Genuinely
extra-conceptualist spaces like the dual of $L^\infty[0,1]$ are seen as
pathological and do not attract substantial mainstream interest.

The naive concept of a set is not only grounded in an untenable philosophical
stance, it also yields a picture of a universe which bears little relation to
mainstream mathematics. When we pay closer attention to the kinds of
collections that we actually use, we find that they stratify
into three levels, and this stratification fits with mathematical
practice in a way that Cantorian set theory does not. For more
on this point see \cite{W5}.

%##

\end{document}